\documentstyle[12pt]{article}
\def\e{{\varepsilon}}

\def\arr{{\longrightarrow}}

\def\q{{\quad$\Box$}}

\def\binfty{{\lefteqn{\hspace*{-0.04em}\infty}{\lefteqn{\hspace*{0.04em}\infty}{\lefteqn{\hspace*{-0.02em}\infty}{\lefteqn{\hspace*{0.02em}\infty}{\infty}}}}}}

\def\be{\begin{equation}}
\def\ee{\end{equation}}

\newtheorem{thm}{Theorem}[section]
\newtheorem{lem}[thm]{Lemma}
\newtheorem{prop}[thm]{Proposition}
\newtheorem{dfn}[thm]{Definition}
\newtheorem{cor}[thm]{Corollary}

\newcommand{\norm}[1]{\left\| #1 \right\|}

\date{}
\author{V.~M.~Manuilov}
\title{On $C^*$-algebras related to asymptotic homomorphisms
\footnote{This research was partially supported by
RFBR (grant No 99-01-01201) and by INTAS (grant No 96-1099).}}

\begin{document}

\maketitle

 \begin{abstract}
We study the $C^*$-algebras related to Mishchenko's version of asymptotic
homomorphisms. In particular we show that their different versions are
weakly homotopy equivalent but not isomorphic to each other. We give also
the continuous version for these algebras.

 \end{abstract}

{\protect\small\protect
\section{Introduction}

}

Let $M_k$, $k\in{\bf N}$, denote the $k{\times}k$
matrix $C^*$-algebra.
To make notation uniform we will write $M_\infty$ for the $C^*$-algebra $K$
of compact operators. Let $\kappa=\{\kappa(1),\kappa(2),\ldots\}$ be a
function from ${\bf N}$ to ${\bf N}\cup\{\infty\}$, i.e. $\kappa$ is a
sequence with $\kappa(i)$, $i\in\bf N$, being either integers or infinity.
We suppose that the sequences $\kappa=\{\kappa(i)\}$ are {\it monotonely
non-decreasing\/} and either stabilize (there exists
$\lim_{i\to\infty}\kappa(i)$) or approach infinity as $i\to\infty$.

Consider the quotient $C^*$-algebra
 $$
\prod_{i=1}^\infty M_{\kappa(i)}/\oplus_{i=1}^\infty M_{\kappa(i)},
 $$
which is the target $C^*$-algebra for discrete asymptotic homomorphisms
\cite{Connes-Higson,Loring},
i.e. any discrete asymptotic homomorphism from a $C^*$-algebra $A$
into matrix algebras or into compacts gives a genuine $*$-homomorphism
from $A$ into this target algebra.
In \cite{mish-noor} Mishchenko suggested to consider another version of
discrete asymptotic homomorphisms, namely to add the following property: if
$\varphi_n:A\arr B$ is a discrete asymptotic homomorphism, it should satisfy
also
 $$
\lim_{i\to\infty}\norm{\varphi_{i+1}(a)-\varphi_i(a)}=0
 $$
for any $a\in A$.
The corresponding target $C^*$-algebras were constructed in \cite{man-mish}
as follows. Since any sequence $\kappa$ is non-decreasing, there is a natural
inclusion $M_{\kappa(i)}\subset M_{\kappa(i+1)}$, hence the right shift
\begin{equation}\label{alpha}
\alpha:(m_1,m_2,\ldots)\longmapsto(0,m_1,m_2,\ldots),
\end{equation}
where $m=(m_1,m_2,\ldots)\in \prod_{i=1}^\infty M_{\nu(i)}$,
is well-defined both on $\prod_{i=1}^\infty M_{\kappa(i)}$ and on
$\prod_{i=1}^\infty M_{\kappa(i)}/\oplus_{i=1}^\infty M_{\kappa(i)}$. The
map $\alpha$ is an endomorphism of these $C^*$-algebras, and it is easy to
see that the $\alpha$-invariant subset in $\prod_{i=1}^\infty
M_{\kappa(i)}/\oplus_{i=1}^\infty M_{\kappa(i)}$ is a $C^*$-algebra too.
We denote this $\alpha$-invariant $C^*$-algebra by $Q(\kappa)$.

If the sequence $\kappa$ stabilizes, $\lim_{i\to\infty}\kappa(i)=n$,
$n\in{\bf N}\cup\{\infty\}$, we denote such $\kappa$ by $\bf n$ (or by
$\binfty$ if $n=\infty$). In this case we write $Q({\bf n})$ instead of
$Q(\kappa)$. The $C^*$-algebra $Q(\kappa)$ can be described as a set
$B(\kappa)$ of all sequences $(m_1,m_2,\ldots)$ of
matrices $m_i\in M_{\kappa(i)}$ such that
 $$
\lim_{i\to\infty}\norm{m_{i+1}-m_i}=0
 $$
modulo the sequences converging to zero. The set $B(\kappa)$ is a
$C^*$-algebra too, so $Q(\kappa)$ is a quotient $C^*$-algebra of
$B(\kappa)$ modulo the ideal $I(\kappa)$ of sequences converging to zero.
Moreover, the ideal $I(\kappa)$ is an essential ideal in $B(\kappa)$, so
$B(\kappa)$ lies in the multiplier $C^*$-algebra $M(I(\kappa))$.

Note that though the $C^*$-algebras $Q(\kappa)$ consist only of
$\alpha$-invariant elements of $\prod_i M_{\kappa(i)}/\oplus_i
M_{\kappa(i)}$, one should not think that $Q(\kappa)$ are smaller than
the latter algebras.
For example, consider $Q(\binfty)$ and let $r=r(n)$ be a sequence of
increasing integers with
\begin{equation}\label{infinit}
\lim_{i\to\infty}(r(i+1)-r(i))=\infty
\end{equation}
then it is easy to see that the sequences $(m_1,m_2,\ldots)$ with
$m_{r(i)}=0$ form an ideal $I_r$ of $Q(\binfty)$ and the
corresponding quotient $C^*$-algebra $Q(\binfty)/I_r$ can be mapped to
$\prod_i M_\infty/\oplus_i M_\infty$ by sending $(m_1,m_2,\ldots)$ to
$(m_{r(1)},m_{r(2)},\ldots)$. Using (\ref{infinit}) one can lift
any element of $\prod_i M_\infty/\oplus_i M_\infty$ by linear
interpolation, hence this map is epimorphic.

Remark that the $C^*$-algebra $Q({\bf 1})$ is commutative, and its spectrum
is the Higson corona, i.e. the quotient space of the Stone--\v{C}ech corona
of $\bf N$ modulo the standard action of $\bf Z$. It is easy to see that the
$C^*$-algebra $Q(\kappa)$ is unital if and only if $\kappa={\bf n}$ for
some finite $n$. In this case one has $Q({\bf n})=Q({\bf 1})\otimes M_n$.
Otherwise $Q(\kappa)$ contains the tensor product $Q({\bf 1})\otimes
M_\infty$ as a proper subalgebra. Indeed, the lifting of $Q({\bf
1})\otimes M_\infty$ to $B(\kappa)$ consists of uniformly vanishing
sequences $m=(m_1,m_2,\ldots)\in B(\kappa)$. This means that for any
such $m$ and for any $\e>0$ there exists a number $j$ such that
 $$
\norm{(m_1,m_2,\ldots)-
(p_jm_1p_j\ ,\ p_jm_2p_j\ ,\ \ldots)}<\e,
 $$
where $p_j\in M_\infty$ denotes a projection onto the first $j$
coordinates. For any $\kappa$ approaching
infinity it is easy to find an element of $B(\kappa)$, which is not
uniformly vanishing. To do so consider a sequence $\{n_i\}$ with
$\lim_{i\to\infty}(n_{i+1}-n_i)=\infty$ and a non-zero element
$\lambda=(\lambda_1,\lambda_2,\ldots)\in B({\bf 1})$ such that
$\lambda_{n_i}=0$ for all $i\in\bf N$. Let $\chi_{[n_i,n_{i+1}]}$ be the
characteristic function for the set $[n_i,n_{i+1}]\subset\bf N$. Let
$\{k_i\}$ be another sequence approaching infinity and let $e_i\in
M_\infty$ be the projection onto the $i$-th coordinate. Then the element
$\oplus_i \chi_{[n_i,n_{i+1}]}\lambda e_{k_i}$ does not lie in $B({\bf
1})\otimes M_\infty$.


{\protect\small\protect
\section{Weak homotopic equivalence of $Q(\kappa)$}

}

The $C^*$-algebras $Q(\kappa)$ are quotient algebras modulo an essential
ideal, so it is natural to consider the corresponding weaker notion of
homotopy for $*$-homomorphisms into corona algebras.
Remember that a function $g=g(t):[0,1]\arr B(\kappa)$ is called
continuous with respect to the strict topology \cite{Pedersen} if the
functions $g(t)k(t)$ and $k(t)g(t)$ are continuous in $C([0,1],I(\kappa))$
for any continuous function $k=k(t)\in C([0,1],I(\kappa))$.
The $C^*$-algebra of all continuous functions with respect to the strict
topology we denote by $C_{str}([0,1],B(\kappa))$. The quotient maps
$C_{str}([0,1],B(\kappa))\arr C_{str}([0,1],B(\kappa))/C([0,1],I(\kappa))$
and $B(\kappa)\arr Q(\kappa)$ we denote by $q$.

Let $f_0,f_1:A\arr
Q(\kappa)$ be two $*$-homomorphisms. They are called homotopic
(as maps into quotient algebra) if there
exists a family of maps (not necessarily $*$-homomorphisms) $F_t:A\arr
B(\kappa)$ such that $F_t(a)$ is continuous with respect to the strict
topology for any $a\in A$,
the composition $q\circ F_t:A\arr C([0,1],Q(\kappa))$ is a
$*$-homomorphism and $q\circ F_i$ coincides with $f_i$, $i=0,1$.

As we are dealing with non-separable $C^*$-algebras, we have to consider
even weaker homotopy equivalence. Let $A$, $B$ be $C^*$-algebras. By
$[A,B]$ we denote the set of all $*$-homomorphisms from $A$ to $B$.
\begin{dfn}
{\rm
We call $A$ and $B$ {\it weakly homotopy equivalent\/} if for any {\it
separable\/} $C^*$-algebra $C$ there are {\it natural\/} (with respect to
$C$) maps
 $$
\phi=\phi_C:[C,A]\arr[C,B]\quad{\rm and}\quad \psi=\psi_C:[C,B]\arr[C,A]
 $$
such that they preserve homotopy of $*$-homomorphisms and that
for any $f:C\arr A$ and for any $g:C\arr B$ the compositions
$(\psi\circ\phi)(f)$ and $(\phi\circ\psi)(g)$ are homotopic to $f$ and
$g$ respectively.
}
\end{dfn}
\begin{lem}
If $A$ and $B$ are weakly homotopy equivalent separable $C^*$-algebras
then they are homotopy equivalent.
\end{lem}
{\bf Proof.}
Put $f=\phi_A(id_A):A\arr B$ and $g=\psi_B(id_B):B\arr A$. By $f^*$ we
denote the map $[B,X]\arr[A,X]$ generated by the $*$-homomorphism $f$.
Then one has $f=f^*(id_B)$ and naturality of $\psi$ implies that
 $$
\psi_A\phi_A(id_A)=\psi_A\ f^*(id_B)=f^*\psi_B(id_B)=f^*g=g\circ f
 $$
so by definition $g\circ f$ is homotopic to
$id_A$. By the same way one can show that $f\circ g$ is homotopic to
$id_B$. \q

Though this weak homotopy equivalence is very weak indeed, it preserves
the $K$-theory.
\begin{lem}\label{K-groups}
If $C^*$-algebras $A$ and $B$ are weakly homotopy equivalent then
$K_*(A)\cong K_*(B)$.
\end{lem}
{\bf Proof.}
It is easy to see that weak homotopy equivalence of $C^*$-algebras $A$ and
$B$ implies the same equivalence of $C^*$-algebras $M_\infty\otimes A$ and
$M_\infty\otimes B$. The groups $K_0$ and $K_1$ are defined by homotopy
classes of $*$-homomorphisms into these algebras of the separable
$C^*$-algebras $\bf C$ and $C_0({\bf R})$ respectively. \q

The $K$-groups of the $C^*$-algebras $Q(\kappa)$ were calculated in
\cite{man-mish}.
\begin{prop}[\cite{man-mish}]
For any $\kappa$ one has $K_0(Q(\kappa))={\bf Z}$ and $K_1(Q(\kappa))=0$.
$K_0(Q(\kappa))$ is generated by the sequence $(e_1,e_1,e_1,\ldots)$,
where $e_1$ is the projection onto the first coordinate. \q
\end{prop}
\begin{prop}\label{homotequiv}
Let $\kappa_1$ and $\kappa_2$ be two sequences approaching infinity but
$\kappa_j\neq\binfty$, $j=1,2$. Then the $C^*$-algebras $Q(\kappa_1)$ and
$Q(\kappa_2)$ are homotopy equivalent.
\end{prop}
{\bf Proof.}
Let $\kappa_3$ be a sequence majorizing both $\kappa_1$ and $\kappa_2$.
Transitivity of homotopy equivalence shows that we can prove
homotopy equivalence only for $Q(\kappa_1)$ and $Q(\kappa_3)$. As
$\kappa_1(i)\leq\kappa_3(i)$ for any $i\in\bf N$, there exists a natural
inclusion $f:Q(\kappa_1)\arr Q(\kappa_3)$. To define a map $g$ in the other
direction consider an element $m=(m_1,m_2,\ldots)\in B(\kappa_3)$,
define numbers $r=r(i)$ by $\kappa_1(i-1)<\kappa_3(r)\leq\kappa_1(i)$ and put
$(g(m))_i=m_r$. Then $g(m)\in B(\kappa_1)$ and $g(m)\in I(\kappa_1)$
whenever $m\in I(\kappa_3)$, so the map $g:Q(\kappa_3)\arr
Q(\kappa_1)$ is well-defined. Then the composition $g\circ
f:Q(\kappa_1)\arr Q(\kappa_1)$ is reduced to renumbering: a sequence
$m=(m_1,m_2,\ldots)\in B(\kappa_1)$ is mapped to a sequence
\begin{equation}\label{seq1}
(\underbrace{m_1\ ,\ \ldots\ ,\ m_1}_{s_1\ {\rm times}}\ ,\
\underbrace{m_2\ ,\ \ldots\ ,\ m_2}_{s_2\ {\rm times}}\ ,\ \ldots\ ).
\end{equation}
We are going to construct the homotopy $F_t:Q(\kappa_1)\arr
C_{str}([0,1],B(\kappa_1))$ connecting $g\circ f$ with $id_{Q(\kappa_1)}$.
For any $a\in Q(\kappa_1)$ choose a
representative $m\in B(\kappa_1)$ of the form (\ref{seq1}).
For $t\in[0,1/2]$ define $F_t(a)$ by connecting (by a
linear path) the sequence (\ref{seq1}) with the sequence
\begin{equation}\label{seq2}
F_{1/2}(a)=(m_1\ ,\ \underbrace{m_2\ ,\ \ldots\ ,\ m_2}_{s_1+s_2-1\ {\rm
times}}\ ,\ \underbrace{m_3\ ,\ \ldots\ ,\ m_3}_{s_3\ {\rm
times}}\ ,\ \ldots\ ).
\end{equation}
For $t\in[1/2,3/4]$ we connect the sequence (\ref{seq2}) with the sequence
 $$
F_{3/4}(a)=(m_1\ ,\ m_2\ ,\ \underbrace{m_3\ ,\
\ldots\ ,\ m_3}_{s_1+s_2+s_3-1\ {\rm times}}\ ,\
\underbrace{m_4\ ,\ \ldots\ ,\ m_4}_{s_4\ {\rm times}}\ ,\ \ldots\ ),
 $$
etc. Finally we put $F_1(a)=(m_1,m_2,\ldots)$.
It is easy to see that $F_t(a)$ is continuous with respect to the strict
topology (one has to check this only at $t=1$), that $q\circ F_t(a)$
does not depend on choice of a representative $m$, and that $F_t(\cdot)$
is a $*$-homomorphism modulo $C([0,1],I(\kappa_1))$, so it gives the
necessary homotopy. A similar homotopy connects $f\circ g$ with
$id_{Q(\kappa_3)}$. \q
\begin{prop}\label{konechn}
Let $A$ be a separable $C^*$-algebra and let $f:A\arr Q(\binfty)$ be a
$*$-ho\-mo\-mor\-phism. Then there exists a {\em finite} sequence $\kappa$
such that $f(A)\subset Q(\kappa)\subset Q(\binfty)$.
\end{prop}
{\bf Proof.}
Fix a dense sequence $\{a_i\}$, $i\in{\bf N}$, in $A$ and a
monotonely decreasing sequence of numbers $\{\e_n\}$, $n\in\bf N$,
such that $\lim_{n\to\infty}\e_n=0$. Let
$m^{(i)}=\left(m^{(i)}_1,m^{(i)}_2,\ldots\right)\in B(\binfty)$ be
liftings for $f(a_i)$. Then the numbers $\kappa(n)$ can be chosen to
satisfy
\begin{equation}\label{posl}
\norm{m^{(i)}_j-p_{\kappa(n)}m^{(i)}_jp_{\kappa(n)}}<\e_n
\end{equation}
for all $1\leq i,j\leq n$, where $p_{\kappa(n)}\in M_\infty$ are projections
onto the first $\kappa(n)$ coordinates.
It follows from (\ref{posl}) that for any $i$ the sequences
 $$
\left(m^{(i)}_1,\ m^{(i)}_2,\ \ldots\ \right)\quad{\rm and}\quad
\left(p_{\kappa(1)}m^{(i)}_1p_{\kappa(1)}\ ,\
p_{\kappa(2)}m^{(i)}_2p_{\kappa(2)}\ ,\ \ldots\ \right)
 $$
give the same element in
$Q(\binfty)$ and the last sequence lies in $B(\kappa)$. As
$Q(\kappa)\subset Q(\binfty)$ is a $C^*$-subalgebra, $f$ is
continuous, and $f(a_i)\subset Q(\kappa)$, so $f(A)\in Q(\kappa)$ too. \q

Propositions \ref{homotequiv} and \ref{konechn} imply
\begin{thm}
If $\kappa$ approaches infinity then all $C^*$-algebras $Q(\kappa)$ are
weakly homotopy equivalent to each other. \q
\end{thm}

{\protect\small\protect
\section{Growth functions and $Q(\kappa)$}

}

Now we are going to show that the $C^*$-algebras $Q(\kappa)$ can be
different for different growth of $\kappa$.
We say that $\kappa$ has polynomial growth if there exist numbers
$0<C_1<C_2$ and $n\geq 1$ such that
 $$C_1 i^n\leq \kappa(i)\leq C_2 i^n\quad{\rm and}\quad
\lim_i\frac{\kappa(i+1)-\kappa(i)}{i^n}=0.
 $$
It has exponential growth if there exists a number $a>1$ such that for all
$i$ one has $\kappa(i+1)\geq a\,\kappa(i)$.
\begin{thm}\label{nonisom}
Let $\kappa_1$ have polynomial growth and let $\kappa_2$ have exponential
growth.

1.
There exists a non-trivial $Q({\bf 1})$-valued trace on $Q(\kappa_1)$.

2.
There is no non-trivial trace on $Q(\kappa_2)$.
\end{thm}
{\bf Proof.}
1.
We define a trace on the $C^*$-algebra $B(\kappa_1)$ and then show that
it vanishes on the ideal $I(\kappa_1)$.
Let $m=(m_1,m_2,\ldots)\in B(\kappa_1)$, $m_i\in M_{\kappa_1(i)}$.
Put
 $$
\tau(m)=\left(\frac{{\rm tr}_{\kappa_1(1)}(m_1)}{\kappa_1(1)}\ ,\
\frac{{\rm tr}_{\kappa_1(2)}(m_2)}{\kappa_1(2)}\ ,\ \frac{{\rm
tr}_{\kappa_1(3)}(m_3)}{\kappa_1(3)}\ ,\ \ldots\ \right),
 $$
where ${\rm tr}_i$ is the standard trace on $M_i$ normalized by ${\rm
tr}_i(1_i)=i$ for $1_i$ being the identity of $M_i$.
If $\norm{m}=1$ then
${\rm tr}_{\kappa_1(i)}(m_i)\leq \kappa_1(i)$, so $\norm{\tau(m)}\leq 1$.
One has also
\begin{eqnarray*}
\norm{\frac{{\rm tr}_{\kappa_1(i+1)}(m_{i+1})}{\kappa_1(i+1)}-\frac{{\rm
tr}_{\kappa_1(i)}(m_i)}{\kappa_1(i)}}&\leq& \norm{\frac{{\rm
tr}_{\kappa_1(i+1)}(m_{i+1})-{\rm tr}_{\kappa_1(i)}(m_i)}{\kappa_1(i+1)}}+
\norm{\frac{{\rm tr}_{\kappa_1(i)}(m_i)}{\kappa_1(i)\kappa_1(i+1)}}\\
&\leq&\norm{m_{i+1}-m_i}+\frac{1}{\kappa_1(i+1)}\to 0,
\end{eqnarray*}
hence the image of $\tau$ lies in $B({\bf 1})$. Finally the property
$\tau(ab)=\tau(ba)$ for $a,b\in B(\kappa_1)$ is obvious.
So $\tau$ is a $B({\bf 1})$-valued trace on $B(\kappa_1)$.
It is also obvious that $\tau(I(\kappa_1))\subset I({\bf 1})$, so $\tau$ is
well-defined on the quotient algebra $Q(\kappa_1)$. It remains to show
that $\tau$ is non-trivial. Let
 $$
m_i={\rm diag}\left\lbrace
1\ ,\ \frac{\kappa_1(i)-1}{\kappa_1(i)}\
,\ \ldots\ ,\ \frac{1}{\kappa_1(i)}\right\rbrace
 $$
be diagonal matrices in $M_{\kappa_1(i)}$.
Then
 $$
\norm{m_{i+1}-m_i}=\frac{\kappa_1(i+1)-\kappa_1(i)}{\kappa_1(i+1)}
\leq \frac{1}{C_1}\frac{\kappa_1(i+1)-\kappa_1(i)}{i^n}\to 0,
 $$
so the sequence $m=(m_1,m_2,\ldots)$ lies in $B(\kappa_1)$, hence it
defines an element in $Q(\kappa_1)$ and one has
 $$
\tau(m)=\left(\frac{\kappa_1(1)+1}{2\kappa_1(1)}\ ,\
\ldots\ ,\ \frac{\kappa_1(i)+1}{2\kappa_1(i)}\ ,\ \ldots\ \right).
 $$
As $\lim_i\frac{\kappa_1(i)+1}{2\kappa_1(i)}=\frac{1}{2}\neq 0$, the trace
$\tau$ is non-trivial.

2.
For convenience of notation we assume that $a=2$ and that
$\kappa_2(i+1)\geq 2\kappa_2(i)$.
We also assume that all $M_i$ are embedded in the
$C^*$-algebra of bounded operators on a Hilbert space $H$ with a fixed basis
$\xi_1,\xi_2,\ldots$. Let $u_1$, $u_2$ be isometries on $H$ such that
$u_1(\xi_k)=\xi_{2k-1}$, $u_2(\xi_k)=\xi_{2k}$, $k\in\bf N$.
Let $m=(m_1,m_2,\ldots)\in B(\kappa_2)$. Then one can define an element
$m\oplus m\in B(\kappa_2)$ as follows. Put $(m\oplus m)_1=0$ and
 $$
(m\oplus m)_{i+1}=u_1 m_i u_1^*+u_2 m_i u_2^*
 $$
for $i>1$. Then $\dim (m\oplus m)_{i+1}=2\kappa_2(i)\leq\kappa_2(i+1)$,
and it is easy to see that $\norm{(m\oplus m)_{i+1}-(m\oplus m)_i}\to 0$.
As $\norm{(m\oplus m)_i}\to 0$ whenever $m\in I(\kappa_2)$, so the map
$m\mapsto (m\oplus m)$ defines a homomorphism $Q(\kappa_2)\arr Q(\kappa_2)$.
In a similar way one can define a direct sum of 4, 8, etc. summands.

Now suppose that there exists a trace $\tau$ on $Q(\kappa_2)$ with the
property $\tau(ab)=\tau(ba)$ for $a,b\in Q(\kappa_2)$. Then $\tau$ extends
to a trace on $B(\kappa_2)$ and we still denote this extension by $\tau$.
Let $m\in B(\kappa_2)$ be a positive element with $\norm{m}=1$ and
$\tau(m)\neq 0$.
One has $m\oplus m=m^{(1)}+m^{(2)}$, where
 $$
m^{(1)}=(0\ ,\ u_1 m_1 u_1^*\ ,\ u_1 m_2 u_1^*\ ,\ u_1 m_3 u_1^*\ ,\
\ldots\ ),
 $$
 $$
m^{(2)}=(0\ ,\ u_2 m_1 u_2^*\ ,\ u_2 m_2 u_2^*\ ,\ u_2 m_3 u_2^*\ ,\
\ldots\ ).
 $$
As all $m_i$ are positive and of dimension $\kappa(i)$, so the sequences
 $$
a_j=(0\ ,\ u_j (m_1)^{1/2}\ ,\ u_j (m_2)^{1/2}\ ,\ u_j (m_3)^{1/2}\ ,\
\ldots\ )
 $$
and
 $$
b_j=(0\ ,\ (m_1)^{1/2} u_j^*\ ,\ (m_2)^{1/2} u_j^*\ ,\ (m_3)^{1/2} u_j^*\ ,\
\ldots\ ),
 $$
$j=1,2$, are well-defined as elements of the algebra $B(\kappa_2)$.
One has
 $$
\tau(m^{(j)})=\tau(a_jb_j)=\tau(b_ja_j)=\tau(\alpha(m)),
 $$
where $\alpha$ is the shift (\ref{alpha}), so $\tau(m^{(j)})=\tau(m)$, as
$m=\alpha(m)$ modulo $I(\kappa_2)$. Therefore $\tau(m\oplus m)=2\tau(m)$.
By the same way we can show that for any $k$ one has
 $$
\underbrace{m\oplus m\oplus \ldots\oplus m}_{2^k\ {\rm times}}
\in B(\kappa_2)
 $$
and
 \begin{equation}\label{estim1}
\tau(\underbrace{m\oplus m\oplus \ldots\oplus m}_{2^k\ {\rm
times}})=2^k\tau(m).
 \end{equation}
On the other hand it is easy to see that
 $$
\|\underbrace{m\oplus m\oplus \ldots\oplus m}_{2^k\ {\rm times}}\|=
\norm{m},
 $$
hence
 \begin{equation}\label{estim2}
\tau(\underbrace{m\oplus m\oplus \ldots\oplus m}_{2^k\
{\rm times}})\leq 1
 \end{equation}
and the contradiction between (\ref{estim1}) and (\ref{estim2})
shows that $\tau(m)=0$ and finishes the proof. \q

Remember that a $C^*$-algebra $A$ is called stable if $A\cong A\otimes
M_\infty$.
\begin{cor}
If $\kappa$ is of polynomial growth $($or smaller$)$ then the $C^*$-algebra
$Q(\kappa)$ is not stable. If $\kappa$ is of exponential growth
then $Q(\kappa)$ is stable. \q
\end{cor}
\begin{cor}\label{nol}
Let $\kappa$ be of polynomial growth. Let $\tau$ be a trace constructed in
Theorem \ref{nonisom} and let $a\in Q({\bf 1})\otimes
M_\infty\subset Q(\kappa)$. Then $\tau(a)=0$. \q
\end{cor}
Let $\kappa$ be of polynomial growth. Consider a $C^*$-subalgebra
$I_\tau\subset Q(\kappa)$ generated by all positive elements $a\in
Q(\kappa)$ with $\tau(a)=0$. As
 $$
\tau(ba)=\tau(a^{1/2}ba^{1/2})\leq\norm{b}\tau(a)=0
 $$
for any $b\in Q(\kappa)$, so $I\tau$ is an ideal in $Q(\kappa)$.
By Lemmas \ref{K-groups} and \ref{nol} the $K$-groups of $I_\tau$ coincide
with those of $Q(\kappa)$ and the inclusion $I_\tau\subset Q(\kappa)$
induces an isomorphism of their $K$-groups. Hence the six-term exact
sequence shows that the quotient $C^*$-algebra $Q(\kappa)/I_\tau$ is
$K$-contractible.

{\protect\small\protect
\section{Continuous version for $Q(\kappa)$}

}

Let $C_b([0,\infty),B)$ (resp. $C_0([0,\infty),B)$) be the $C^*$-algebra
of bounded continuous (resp. vanishing at infinity) functions on the
half-line taking values in a $C^*$-algebra $B$. The quotient algebra
$Q_b(B)=C_b([0,\infty),B)/C_0([0,\infty),B)$ is usually used for describing
continuous asymptotic homomorphisms. But this algebra is too big and one
has often to make reparametrization in order to obtain slowly varying
functions. Consider at first the case $B=M_\infty$. There exists a natural
inclusion
 \begin{equation}\label{otobrfactoralgebr}
Q(\binfty)\arr Q_b(M_\infty),
 \end{equation}
given by the natural inclusion ${\bf N}\subset [0,\infty)$
and by linear interpolation from sequences of compact operators to the
$M_\infty$-valued functions. This inclusion is obviously not epimorphic.

For a function $f(t)\in C_b([0,\infty),M_\infty)$ and for an interval
$[a,b]\subset[0,\infty)$ put
 $$
{\rm Var}_{a,b}f=\sup_{t_1,t_2\in[a,b]}\norm{f(t_2)-f(t_1)}
 $$
and consider in $C_b([0,\infty),M_\infty)$ a subset of functions with
variation vanishing at infinity
 $$
C_b^{vv}([0,\infty),M_\infty) =\{f(t)\in C_b([0,\infty),M_\infty):
\lim_{t\to\infty}{\rm Var}_{a+t,b+t}f=0\}.
 $$
This subset is a $C^*$-subalgebra and it does not depend on choice of $a$
and $b$. Denote by $Q_b^{vv}(M_\infty)$ the quotient $C^*$-algebra
$C_b^{vv}([0,\infty),M_\infty)/C_0([0,\infty),M_\infty)$.
 \begin{lem}\label{isom-alg}
$C^*$-algebras $Q(\binfty)$ and $Q_b^{vv}(M_\infty)$ are isomorphic.
 \end{lem}
{\bf Proof.}
The image of $Q(\binfty)$ under the map (\ref{otobrfactoralgebr})
obviously lies in $Q_b^{vv}(M_\infty)$.
To get the map in the other direction one has to assign to a function
$f(t)\in C_b^{vv}([0,\infty),M_\infty)$ representing an element of
$Q_b^{vv}(M_\infty)$ a sequence $\{f(n)\}$ of its values at integer points.
One can directly check that these $*$-homomorphisms are inverse to each
other. \q
\begin{prop}
The $C^*$-algebras $Q_b(M_\infty)$ and $Q_b^{vv}(M_\infty)$ are weakly
homotopy equivalent.
\end{prop}
{\bf Proof.}
It is easy to see that for any $*$-homomorphism of a separable
$C^*$-algebra $A$ into $Q_b(M_\infty)$ one can find a
reparametrization $s=s(t)$ such that after passing from $t$ to $s$ the
image of $A$ lies in $Q_b^{vv}(M_\infty)$. \q

To imitate the continuous version for general $\kappa$ one should consider
telescopes. Let
 $$
T(\kappa)=\{f(t)\in C_b([0,\infty),M_\infty):
f(t)\in M_{\kappa(i)}\ \ {\rm for}\ t\in [i,i+1)\}
 $$
be the telescope $C^*$-algebra. Put
 $$
C_0([0,\infty),\kappa)
=\{f(t)\in T(\kappa):
\lim_{t\to\infty}\norm{f(t)}=0\}
 $$
and
 $$
C_b^{vv}([0,\infty),\kappa)
=\{f(t)\in T(\kappa):
\lim_{t\to\infty}{\rm Var}_{t,t+1}f=0\}.
 $$
 \begin{lem}
$C^*$-algebras $Q(\kappa)$ and
$C_b^{vv}([0,\infty),\kappa)/C_0([0,\infty),\kappa)$ are isomorphic. \q
 \end{lem}


\vspace{2cm}
\noindent
V.~M.~Manuilov\\
Dept. of Mech. and Math.,\\
Moscow State University,\\
Moscow, 119899, RUSSIA\\
e-mail: manuilov@mech.math.msu.su


\begin{thebibliography}{9}

{\small

\bibitem{Connes-Higson}
{\sc A.~Connes, N.~Higson}. Deformations, morphismes asymptotiques et
$K$--theorie bivariante.
{\it C. R. Acad. Sci. Paris}, s\'erie I, {\bf 311} (1990),
101--106.

\bibitem{Loring}
{\sc T.~A.~Loring}. Almost multiplicative maps between $C^*$-algebras.
Operator Algebras and Quantum Field Theory. Internat. Press, 1997,
111--122.


\bibitem{man-mish}
{\sc V. M. Manuilov, A. S. Mishchenko}. Asymptotic and Fredholm
representations of discrete groups. {\it Mat. Sb.\/} {\bf 189} (1998), No 10,
53--72 (in Russian). English translation: {\it Russian Acad. Sci. Sb. Math.\/}
    {\bf 189} (1998), 1485--1504.


\bibitem{mish-noor}
{\sc A.~S.~Mishchenko, Noor Mohammad}. Asymptotic representations of
discrete groups. ``Lie Groups and Lie Algebras. Their Representations,
Generalizations and Applications.'' Mathematics and its Applications, {\bf
433}. Kluver Acad. Publ. Dordrecht, 1998, 299--312.

\bibitem{Pedersen}
{\sc G.~K.~Pedersen.}
$C^*$-algebras and their automorphism groups. ---
London--New York--San Francisco: Academic Press, 1979.




}
\end{thebibliography}
\end{document}